\documentclass[12pt,letterpaper]{amsart}
\usepackage{epic,eepic,latexsym, amssymb, amscd, amsfonts, xypic, floatflt}

 \newlength{\baseunit}               
 \newcount{\numlines}                
\newcommand{\point}{\vspace{3mm}\par \noindent \refstepcounter{subsection}{\bf \thesubsection.} }

\newcommand{\bpoint}[1]{\vspace{3mm}\par \noindent \refstepcounter{subsection}{\bf \thesubsection.} 
  {\bf #1.} }

\newcommand{\SSS}{S}
\newcommand{\one}{\mathbf{1}}
\newcommand{\aaa}{\mathbf{a}}
\newcommand{\bb}{\mathbf{b}}
\newcommand{\cc}{\mathbf{c}}
\newcommand{\dd}{\mathbf{d}}
\newcommand{\ee}{\mathbf{e}}
\newcommand{\ff}{\mathbf{f}}
\newcommand{\xx}{\mathbf{x}}
\newcommand{\yy}{\mathbf{y}}

\newcommand{\Q}{\mathbb{Q}}

\newcommand{\F}{\mathbb{F}}

\newcommand{\proj}{\mathbb P}
\newcommand{\oh}{{\mathcal{O}}}

\newcommand{\al}{\alpha}

\newcommand{\be}{\beta}
\newcommand{\ep}{\epsilon}
\newcommand{\ga}{\gamma}

\newcommand{\si}{\sigma}

\newcommand{\Segre}{\mathcal{S}_3}
\newcommand{\Igusa}{\mathcal{I}_4}

\newcommand{\cited}{}

\newcommand{\lremind}[1]{{\bf[label:  #1]}}
\newcommand{\secretnote}[1]{}
\newcommand{\notation}[1]{}
\renewcommand{\lremind}[1]{{}}

\newcommand{\cut}[1]{}

\begin{document}
\pagestyle{plain}
\title{{\large {
A description of the outer automorphism of $\SSS_6$, and the
invariants of six points in projective space
}}
}
\author{Ben Howard, John Millson, Andrew Snowden, and Ravi Vakil
}
\email{vakil@math.stanford.edu}
\thanks{Partially supported by NSF CAREER/PECASE Grant DMS--0228011.}
\date{Hallowe'en, 2007.}
\begin{abstract}
  We use a simple description of the outer automorphism of $\SSS_6$ to
  cleanly describe the invariant theory of six points in $\proj^1$,
  $\proj^2$, and $\proj^3$.
\end{abstract}
\maketitle

{\parskip=12pt 
  
  In \S \ref{outer1}--\ref{outer2}, we give two short descriptions of
  the outer automorphism of $\SSS_6$, complete with a proof that they
  indeed describe an outer automorphism.  (Our goal is to have a
  construction that the reader can fully understand and verify while
  sitting on a bus, perhaps doodling in the margins.)  In \S
  \ref{outer3} we give another variation on this theme that is
  attractive, but more long-winded.  The latter two descriptions do
  not distinguish any of the six points.  Of course, these
  descriptions are equivalent to the traditional one (\S \ref{usual})
  --- there is after all only one nontrivial outer automorphism
  (modulo inner automorphisms).  In \S \ref{invariant}, we use this to
  cleanly describe the invariant theory of six points in projective
  space.  This is not just a random application; the descriptions of
  \S \ref{outersection} were discovered by means of this invariant
  theory.  En route we use the outer automorphism to describe
  five-dimensional representations of $\SSS_5$ and $\SSS_6$, \S
  \ref{reptheory}.

The outer automorphism was first described by H\"older in 1895.  Most
verifications use some variation of Sylvester's synthemes, or work
directly with generators of $\SSS_6$; non-trivial calculation is often
necessary.  Other interpretations are in terms of finite geometries,
for example involving finite fields with $2$, $3$, $4$, $5$, or $9$
elements, and are beautiful, but require non-trivial verification.

\section{The outer automorphism of $\SSS_6$}
\label{outersection}

\bpoint{First description of the outer automorphism: the mystic
  pentagons} Consider a complete graph on five vertices numbered $1$
through $5$.  The reader will quickly verify that there are precisely
six ways to two-color the edges (up to choices of colors) so that the
edges of one color (and hence the other color) form a $5$-cycle, see
Figure~\ref{pentagon}.  We dub these the six mystic pentagons.  Then
$\SSS_5$ acts on the six mystic pentagons by permuting the vertices,
giving a map $i: \SSS_5 = \SSS_{\{1,\dots, 5 \}} \rightarrow \SSS_{ \{
  \aaa, \dots, \ff \}} = \SSS_6$.  This is an inclusion --- the kernel
must be one of the normal subgroups  $\{ e \}$, $A_5$, or $\SSS_5$, but we visually 
verify that $(123)$ acts nontrivially.  Moreover, it is not a usual
inclusion as $(12)$ induces permutation $(\aaa \dd) (\bb \cc)
(\ee \ff)$ --- not a transposition.  Hence $\SSS_6 = \SSS_{ \{ \aaa, \dots, \ff \} }$ acts on
the six cosets of $i(\SSS_5)$, inducing a map $f: \SSS_{ \{ \aaa,
  \dots, \ff\} } \rightarrow \SSS_{ \{ 1, \dots, 6 \}}$.  This is the
outer automorphism.  This can be verified in several ways (e.g.,
$(\aaa \dd) (\bb \cc) (\ee \ff)$ induces the nontrivial permutation
$(12) \in \SSS_{ \{ 1, \dots, 6 \}}$, so $f$ is injective and hence an
isomorphism; and $i$ is not a usual inclusion, so $f$ is not inner),
but for the sake of simplicity we do so by way of a second
description.
\label{outer1}\lremind{outer1}

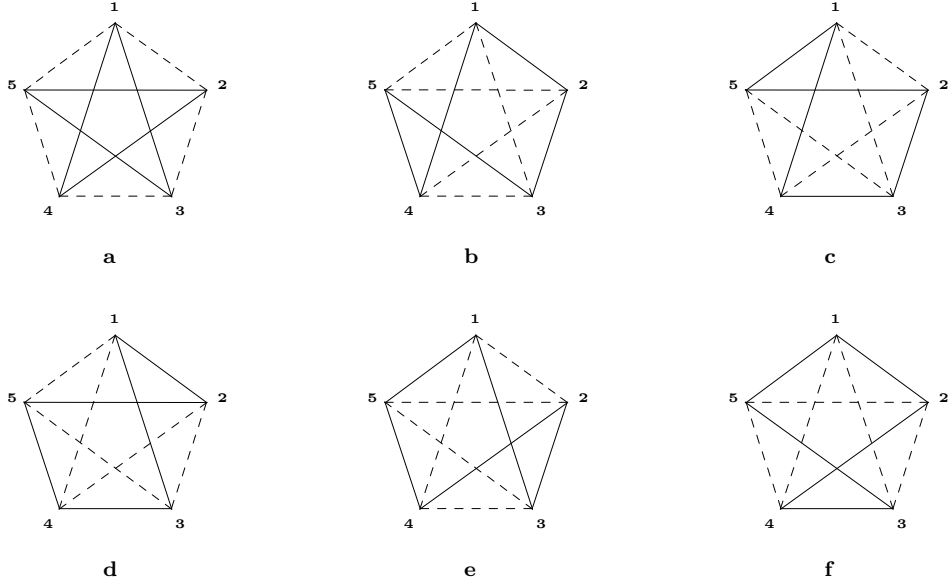
\begin{figure}[ht]
\begin{center}
\setlength{\unitlength}{0.00083333in}
\begingroup\makeatletter\ifx\SetFigFont\undefined%
\gdef\SetFigFont#1#2#3#4#5{%
  \reset@font\fontsize{#1}{#2pt}%
  \fontfamily{#3}\fontseries{#4}\fontshape{#5}%
  \selectfont}%
\fi\endgroup%
{\renewcommand{\dashlinestretch}{30}
\begin{picture}(5813,3655)(0,-10)
\path(104,3064)(1247,3063)
\path(1247,3063)(324,2397)
\path(324,2397)(674,3484)
\path(674,3484)(1027,2397)
\path(1027,2397)(104,3066)
\dashline{60.000}(107,3061)(674,3482)
\dashline{60.000}(1243,3064)(1025,2400)
\dashline{60.000}(1025,2400)(327,2400)
\dashline{60.000}(327,2400)(109,3061)
\dashline{60.000}(674,3482)(1243,3064)
\path(104,1114)(1247,1113)
\dashline{60.000}(1247,1113)(324,447)
\dashline{60.000}(324,447)(674,1534)
\path(674,1534)(1027,447)
\dashline{60.000}(1027,447)(104,1116)
\dashline{60.000}(107,1111)(674,1532)
\path(674,1532)(1243,1114)
\dashline{60.000}(1243,1114)(1025,450)
\path(1025,450)(327,450)
\path(327,450)(109,1111)
\dashline{60.000}(2354,1114)(3497,1113)
\path(3497,1113)(2574,447)
\dashline{60.000}(2574,447)(2924,1534)
\path(2924,1534)(3277,447)
\dashline{60.000}(3277,447)(2354,1116)
\path(2357,1111)(2924,1532)
\dashline{60.000}(2924,1532)(3493,1114)
\path(3493,1114)(3275,450)
\dashline{60.000}(3275,450)(2577,450)
\path(2577,450)(2359,1111)
\dashline{60.000}(2354,3064)(3497,3063)
\dashline{60.000}(3497,3063)(2574,2397)
\path(2574,2397)(2924,3484)
\dashline{60.000}(2924,3484)(3277,2397)
\path(3277,2397)(2354,3066)
\dashline{60.000}(2357,3061)(2924,3482)
\path(2924,3482)(3493,3064)
\path(3493,3064)(3275,2400)
\dashline{60.000}(3275,2400)(2577,2400)
\path(2577,2400)(2359,3061)
\path(4604,3064)(5747,3063)
\dashline{60.000}(5747,3063)(4824,2397)
\path(4824,2397)(5174,3484)
\dashline{60.000}(5174,3484)(5527,2397)
\dashline{60.000}(5527,2397)(4604,3066)
\path(4607,3061)(5174,3482)
\dashline{60.000}(5174,3482)(5743,3064)
\path(5743,3064)(5525,2400)
\path(5525,2400)(4827,2400)
\dashline{60.000}(4827,2400)(4609,3061)
\dashline{60.000}(4604,1114)(5747,1113)
\path(5747,1113)(4824,447)
\dashline{60.000}(4824,447)(5174,1534)
\dashline{60.000}(5174,1534)(5527,447)
\path(5527,447)(4604,1116)
\path(4607,1111)(5174,1532)
\path(5174,1532)(5743,1114)
\dashline{60.000}(5743,1114)(5525,450)
\path(5525,450)(4827,450)
\dashline{60.000}(4827,450)(4609,1111)
\put(0,3068){\makebox(0,0)[lb]{\smash{{{\SetFigFont{5}{6.0}{\rmdefault}{\mddefault}{\updefault}$\mathbf{5}$}}}}}
\put(638,3556){\makebox(0,0)[lb]{\smash{{{\SetFigFont{5}{6.0}{\rmdefault}{\mddefault}{\updefault}$\mathbf{1}$}}}}}
\put(1313,3068){\makebox(0,0)[lb]{\smash{{{\SetFigFont{5}{6.0}{\rmdefault}{\mddefault}{\updefault}$\mathbf{2}$}}}}}
\put(1050,2281){\makebox(0,0)[lb]{\smash{{{\SetFigFont{5}{6.0}{\rmdefault}{\mddefault}{\updefault}$\mathbf{3}$}}}}}
\put(225,2281){\makebox(0,0)[lb]{\smash{{{\SetFigFont{5}{6.0}{\rmdefault}{\mddefault}{\updefault}$\mathbf{4}$}}}}}
\put(600,1981){\makebox(0,0)[lb]{\smash{{{\SetFigFont{8}{9.6}{\rmdefault}{\mddefault}{\updefault}$\mathbf{a}$}}}}}
\put(0,1118){\makebox(0,0)[lb]{\smash{{{\SetFigFont{5}{6.0}{\rmdefault}{\mddefault}{\updefault}$\mathbf{5}$}}}}}
\put(638,1606){\makebox(0,0)[lb]{\smash{{{\SetFigFont{5}{6.0}{\rmdefault}{\mddefault}{\updefault}$\mathbf{1}$}}}}}
\put(1313,1118){\makebox(0,0)[lb]{\smash{{{\SetFigFont{5}{6.0}{\rmdefault}{\mddefault}{\updefault}$\mathbf{2}$}}}}}
\put(1050,331){\makebox(0,0)[lb]{\smash{{{\SetFigFont{5}{6.0}{\rmdefault}{\mddefault}{\updefault}$\mathbf{3}$}}}}}
\put(225,331){\makebox(0,0)[lb]{\smash{{{\SetFigFont{5}{6.0}{\rmdefault}{\mddefault}{\updefault}$\mathbf{4}$}}}}}
\put(600,31){\makebox(0,0)[lb]{\smash{{{\SetFigFont{8}{9.6}{\rmdefault}{\mddefault}{\updefault}$\mathbf{d}$}}}}}
\put(2250,1118){\makebox(0,0)[lb]{\smash{{{\SetFigFont{5}{6.0}{\rmdefault}{\mddefault}{\updefault}$\mathbf{5}$}}}}}
\put(2888,1606){\makebox(0,0)[lb]{\smash{{{\SetFigFont{5}{6.0}{\rmdefault}{\mddefault}{\updefault}$\mathbf{1}$}}}}}
\put(3563,1118){\makebox(0,0)[lb]{\smash{{{\SetFigFont{5}{6.0}{\rmdefault}{\mddefault}{\updefault}$\mathbf{2}$}}}}}
\put(3300,331){\makebox(0,0)[lb]{\smash{{{\SetFigFont{5}{6.0}{\rmdefault}{\mddefault}{\updefault}$\mathbf{3}$}}}}}
\put(2475,331){\makebox(0,0)[lb]{\smash{{{\SetFigFont{5}{6.0}{\rmdefault}{\mddefault}{\updefault}$\mathbf{4}$}}}}}
\put(2850,31){\makebox(0,0)[lb]{\smash{{{\SetFigFont{8}{9.6}{\rmdefault}{\mddefault}{\updefault}$\mathbf{e}$}}}}}
\put(2250,3068){\makebox(0,0)[lb]{\smash{{{\SetFigFont{5}{6.0}{\rmdefault}{\mddefault}{\updefault}$\mathbf{5}$}}}}}
\put(2888,3556){\makebox(0,0)[lb]{\smash{{{\SetFigFont{5}{6.0}{\rmdefault}{\mddefault}{\updefault}$\mathbf{1}$}}}}}
\put(3563,3068){\makebox(0,0)[lb]{\smash{{{\SetFigFont{5}{6.0}{\rmdefault}{\mddefault}{\updefault}$\mathbf{2}$}}}}}
\put(3300,2281){\makebox(0,0)[lb]{\smash{{{\SetFigFont{5}{6.0}{\rmdefault}{\mddefault}{\updefault}$\mathbf{3}$}}}}}
\put(2475,2281){\makebox(0,0)[lb]{\smash{{{\SetFigFont{5}{6.0}{\rmdefault}{\mddefault}{\updefault}$\mathbf{4}$}}}}}
\put(2850,1981){\makebox(0,0)[lb]{\smash{{{\SetFigFont{8}{9.6}{\rmdefault}{\mddefault}{\updefault}$\mathbf{b}$}}}}}
\put(4500,3068){\makebox(0,0)[lb]{\smash{{{\SetFigFont{5}{6.0}{\rmdefault}{\mddefault}{\updefault}$\mathbf{5}$}}}}}
\put(5138,3556){\makebox(0,0)[lb]{\smash{{{\SetFigFont{5}{6.0}{\rmdefault}{\mddefault}{\updefault}$\mathbf{1}$}}}}}
\put(5813,3068){\makebox(0,0)[lb]{\smash{{{\SetFigFont{5}{6.0}{\rmdefault}{\mddefault}{\updefault}$\mathbf{2}$}}}}}
\put(5550,2281){\makebox(0,0)[lb]{\smash{{{\SetFigFont{5}{6.0}{\rmdefault}{\mddefault}{\updefault}$\mathbf{3}$}}}}}
\put(4725,2281){\makebox(0,0)[lb]{\smash{{{\SetFigFont{5}{6.0}{\rmdefault}{\mddefault}{\updefault}$\mathbf{4}$}}}}}
\put(5100,1981){\makebox(0,0)[lb]{\smash{{{\SetFigFont{8}{9.6}{\rmdefault}{\mddefault}{\updefault}$\mathbf{c}$}}}}}
\put(4500,1118){\makebox(0,0)[lb]{\smash{{{\SetFigFont{5}{6.0}{\rmdefault}{\mddefault}{\updefault}$\mathbf{5}$}}}}}
\put(5138,1606){\makebox(0,0)[lb]{\smash{{{\SetFigFont{5}{6.0}{\rmdefault}{\mddefault}{\updefault}$\mathbf{1}$}}}}}
\put(5813,1118){\makebox(0,0)[lb]{\smash{{{\SetFigFont{5}{6.0}{\rmdefault}{\mddefault}{\updefault}$\mathbf{2}$}}}}}
\put(5550,331){\makebox(0,0)[lb]{\smash{{{\SetFigFont{5}{6.0}{\rmdefault}{\mddefault}{\updefault}$\mathbf{3}$}}}}}
\put(4725,331){\makebox(0,0)[lb]{\smash{{{\SetFigFont{5}{6.0}{\rmdefault}{\mddefault}{\updefault}$\mathbf{4}$}}}}}
\put(5100,31){\makebox(0,0)[lb]{\smash{{{\SetFigFont{8}{9.6}{\rmdefault}{\mddefault}{\updefault}$\mathbf{f}$}}}}}
\end{picture}
}
\end{center}
\caption{The six mystic pentagons, with black
and white (dashed) edges \lremind{pentagon, fig=mystic2}}
\label{pentagon}
\end{figure}

\bpoint{Second description of the outer automorphism: labeled
  triangles} We now make this construction more symmetric, not
distinguishing the element $6 \in \{ 1, \dots, 6 \}$.  Consider the
$\binom 6 3 = 20$ triangles on six vertices labeled $\{1, \dots, 6
\}$.  There are six ways of dividing the triangles into two sets of
$10$ so that (i) any two disjoint triangles have opposite colors, and (ii)
every tetrahedron has two triangles of each color.  (The bijection
between these and the mystic pentagons $\aaa$, \dots, $\ff$
is as follows.  The triangle $6AB$
is colored the same as edge $AB$.  The triangle $CDE$ ($6 \neq C,D,E$)
is colored the opposite of the ``complementary'' edge $AB$,
where $\{ A,B \} = \{1, \dots, 5 \} - \{ C,D,E \}$.)  The
$\SSS_6$-action on this set is the outer automorphism of
$\SSS_6$. (Reason: $(12)$ induces a nontrivial permutation
$(\aaa \dd) (\bb \cc) (\ee \ff)$ of the mystic
pentagons, so the induced map $\SSS_6
\rightarrow \SSS_{ \{ \aaa, \dots, \ff \} }  \cong \SSS_6$ 
is injective and hence an isomorphism.  But $(12)$
does not induce a transposition on $\{ \aaa, \dots, \ff \}$, so the
automorphism is not inner.)
This isomorphism $\SSS_{ \{ 1, \dots, 6 \} } \rightarrow
\SSS_{\{ \aaa, \dots, \ff \} }$ is inverse to the isomorphism
$f$ of \S \ref{outer1}.
\label{outer2} \lremind{outer2} 

\bpoint{Another description: labeled icosahedra}
\label{outer3}\lremind{outer3}
Here is another description, which is pleasantly $\SSS_6$-symmetric.
Up to rotations and reflections, there are twelve ways to number the
vertices of an icosahedron $1$ through $6$, such that antipodal
vertices have the same label.  Each icosahedron gives ten triples in
$\{ 1, \dots, 6 \}$, corresponding to the vertices around its faces.
These twelve icosahedra come in six pairs, where two icosahedra are
``opposite'' if they have no triples in common.  (It is entertaining
to note that if an icosahedron is embedded in $\Q(\phi)^3$ 
with vertices at $(\pm 1, \pm \phi, 0)$, $(0, \pm 1, \pm
\phi)$, and $(\pm \phi, 0, \pm 1)$, then conjugation in
$\operatorname{Gal}( \Q(\phi)/\Q)$ sends the icosahedron to its
opposite. Here $\phi$ is the golden section.)  Then $\SSS_6$ acts on
these six pairs, and this is the outer automorphism.  One may show
this via bijections to the descriptions of \S \ref{outer1} and \S
\ref{outer2}.  Each pair of mystic icosahedra corresponds to
two-coloring the triangles in $\{1, \dots, 6 \}$, as in \S
\ref{outer2}.  For the bijection to \S \ref{outer1}, the cyclic order
of the vertices around vertex $6$ gives a mystic pentagon.  (This
provides a hands-on way of understanding the $\SSS_6$-action on the
mystic pentagons.)  This description is related to the explanation of
the outer automorphism by John Baez in \cite{baez}.

\bpoint{Relation to the usual description of the outer automorphism of
  $\SSS_6$} The usual description of the outer automorphism is as
follows (e.g.\ \cite{coxeter}).  A {\em syntheme} is
a matching of the numbers $1$, \dots, $6$, i.e.\ an unordered
partition of $\{ 1, \dots, 6 \}$ into three sets of size two.  A {\em
  pentad} is a set of five synthemes whose whose union is the set of
all $15$ pairs.  Then there are precisely six pentads, and the action of
$\SSS_6$ on this set is via the outer automorphism of $\SSS_6$.  We
explain how to get pentads from the mystic pentagons.  Each
mystic pentagon determines a bijection between the white edges and
the black edges, where edge $AB$ corresponds with edge $CD$ if $AB$
and $CD$ don't share a vertex.  If $E = \{ 1, \dots, 5 \} - \{ A,B,C,D
\}$, then to each such pair we obtain the syntheme $AB/CD/E6$, and
there are clearly five such synthemes, no two of which share an edge,
which hence form a pentad.  For example, mystic pentagon $\mathbf{a}$
yields the pentad
$$
\{ 12/35/56, 23/14/56, 34/25/16,  45/13/26, 15/24/36 \}.
$$
\label{usual}\lremind{usual}Another common description of the outer automorphism relates
directly to Figure~\ref{pentagon}.  We find a subgroup $G<S_5$ of size
$20$; we take the subgroup preserving figure $\aaa$ of
Figure~\ref{pentagon}.  Then $S_5$ acts transitively on the six cosets
of $G$, giving a map $i: S_5 \rightarrow S_6$.  This map is an
inclusion as $(123)$ is not in its kernel.  Then $S_6$ acts
(transitively) on the six cosets of $i(S_5)$, yielding a map $\si: S_6
\rightarrow S_6$.  The image (as it is transitive) has size $>2$,
hence (as $S_6$ has only 3 normal subgroups) the kernel is $e$, hence
$\si$ is an automorphism.  Then it is not inner, as $i(S_5)$ is not
one of the six ``obvious'' $S_5$'s in $S_6$.

\bpoint{Representations of $\SSS_5$ and $\SSS_6$} In
Figure~\ref{pentagon}, the edges are colored black and white so that
each edge appears in each color precisely three times with this
choice.  This has the advantage that any odd permutation in $\SSS_5$
(or $\SSS_6$) permutes the six pentagons and exchanges the
colors.  \label{reptheory}\lremind{reptheory}

The pentagons give a convenient way of understanding the two
$5$-dimensional irreducible representations of $\SSS_5$.  The
permutation representation induced by this $\SSS_5$ action on the
mystic pentagons splits into an irreducible $5$-dimensional
representation $F_5$ and a trivial representation $\one$.
The other irreducible $5$-dimensional $\SSS_5$-representation $F'_5$ is
obtained by tensoring $F_5$ with the sign representation $\ep$, which
can be interpreted as the $\SSS_5$ action on the mystic pentagons
``with sign corresponding to color-swapping''.  

There are four irreducible $5$-dimensional representation of $\SSS_6$.
\label{rep6}\lremind{rep6}One 
is the standard representation (which we here denote $B_5$), obtained
by subtracting the trivial representation $\one$ from the usual
permutation representation.  A second is obtained by tensoring with
the sign representation $\ep$: $B'_5 := B_5 \otimes \ep$.  A third is
analogous to the standard representation, obtained by subtracting the
trivial representation $\one$ from the (outer) permutation
representation of $\SSS_6$ on the six mystic pentagons.  One might
denote this the {\em outer automorphism representation}.  The fourth
$5$-dimensional $\SSS_6$-representation is $O'_5 := O_5 \otimes \ep$.
One might term this the {\em signed outer automorphism
  representation}.  It is clear from the construction that 
$F_5$ and $F'_5$ are obtained by
restriction from $O_5$ and $O'_5$.

\point An alternate description of these representations is as
follows.  There are twelve $5$-cycles on vertices labeled $\{1, \dots,
5 \}$, which come in pairs of ``opposites'', consisting of disjoint
$5$-cycles.  Each mystic pentagon is equivalent to such a pair.  If
$\xx$ is a $5$-cycle, denote its opposite by $\overline{\xx}$.  (The
same construction applies for triangles on $\{ 1, \dots, 6 \}$, \S
\ref{outer2}, or labeled icosahedra, \S \ref{outer3}.)  If we have
twelve variables $Z_{\aaa}$, \dots, $Z_{\ff}$, $Z_{\overline{\aaa}}$,
\dots, $Z_{\overline{\ff}}$, with the conditions $Z_{\overline{\xx}}=
- Z_{\xx}$, the $\SSS_6$-action induces the representation $O'_5
\oplus \ep$.  (Of course, the other three representations can be
described similarly, $O_5$ by imposing $Z_{\overline{\xx}}=Z_{\xx}$
instead and replacing $\ep$ by $\one$, $F'_5$ by considering the
$\SSS_5$-action, and $F_5$ by making both
changes.) \label{importantrep}\lremind{importantrep}

\section{The invariant theory of six points in projective space}

We will now relate the outer automorphism of $\SSS_6$ to the space of
six ordered points in projective space, or more precisely the
geometric invariant theory quotient $(\proj^n)^6 // PGL(n+1)$.  The
algebraic statements in this section may be readily checked by any
computer algebra program such as Maple, so the details are omitted.
They were derived using explicit representation theory of $\SSS_6$,
and again the details are unenlightening and will be omitted.
\label{invariant}\lremind{invariant}

In some sense these quotients generalize the notion of cross-ratio,
the space of four points in $\proj^1$.  Any two generic sets of six
ordered points in $\proj^n$ are projectively equivalent if $n>3$, so
the interesting cases are $n=1, 2, 3$.  All three cases were studied
classically, and were known to behave beautifully.

\bpoint{Six points on $\proj^1$}
\label{p1}\lremind{p1}The space of six points on $\proj^1$ may be
interpreted as a threefold in $\proj^5$ cut out by the equations
\cite[p.~17]{do}
$$
Z_1+ \cdots + Z_6 = Z_1^3 + \cdots + Z_6^3 = 0.
$$
(Aside:  This is one of the many ways in which $6$ points
are special. If $m \neq 6$, the space of $m$ points in
$\proj^1$, $(\proj^1)^m // PGL(2)$, is cut out by quadrics
\cite{hmsv}.)
This is the {\em Segre cubic relation}, and this
moduli space is known as the {\em Segre cubic threefold}, which we  denote
$\Segre$.\lremind{ \cite[Example~11.6]{d}}
There is an obvious $\SSS_6$-action on both $(\proj^1)^6$ and the
variables $Z_1$, \dots, $Z_6$.  One might hope that these actions are
conjugate, which would imply some bijection between the six points and
the six variables.  But remarkably, they are related by the outer
automorphism of $\SSS_6$.

An alternate interpretation of this quotient is as the space of
equilateral hexagons in real 3-space, with edges labeled $1$ through
$6$ cyclically, up to translations and rotations \cite{km}.
Rearranging the order of the edges induces a permutation of the
$Z$-variables via the outer automorphism.

Here is a clue that the outer automorphism is relevant.  The
cross-ratio of a certain four of the six points is given by $[Z_1;
\dots; Z_6] \mapsto - (Z_1+Z_2)/ (Z_3+Z_4)$.  A more symmetric avatar
of the cross-ratio of four points on a line is given by
$$
\xymatrix{
(\proj^1)^4  \ar@{-->}[r] & \proj^2 \\ 
(p_1,p_2,p_3,p_4) \ar@{|->} & 
[  (p_2-p_3)( p_1-p_4); (p_1-p_2) (p_3-p_4); (p_1-p_3) (p_4-p_2)
].}
$$
(Here points of $\proj^1$ are written in projective coordinates for
convenience; more correctly we should write $[u_i,v_i]$ for $p_i$,
where $[u_i, v_i] = [p_i,1]$.)    Note that the $\SSS_4$-symmetry is clear in this
manifestation.  
The image is the line $X+Y+Z=0$ in
$\proj^2$.  The traditional cross-ratio is $-X/Y$.  In this
symmetric manifestation, the cross-ratio of a certain four of the six
points is given by $[Z_1; \dots; Z_6] \mapsto [ Z_1+Z_2; Z_3+Z_4;
Z_5+Z_6]$.  The correspondence of a pair of points with 
a ``syntheme'' of the $Z$-variables is a hint that the outer automorphism is somehow
present.

We now describe the moduli map 
$(\proj^1) \dashrightarrow \Segre$ explicitly.
If the points are $p_1, \dots, p_6$ ($1 \leq i \leq 6$), the moduli
map is given (in terms of the second description of the outer automorphism, \S \ref{outer2}) by
$$
Z_{\mathbf{x}} = \sum_{ \{ A, B, C \} \subset \{ 1, \dots, 6 \}} \pm p_A p_B p_C
$$
where $\mathbf{x} \in \{ \aaa, \dots, \ff \}$, and
the sign is $+1$ if triangle $ABC$ is black, and $-1$ if the triangle
is white.
Note that $\sum Z_{\xx} = 0$.
As an added bonus, we see that  the $\SSS_6$-representation on $H^0(\Segre, \oh(1))$ is the signed
outer automorphism representation $O'_5$
 (see \S \ref{importantrep}).

\bpoint{Six points in $\proj^3$, and the Igusa quartic} The Geometric
Invariant Theory
quotient of six points on $\proj^3$ is the Igusa quartic threefold
$\Igusa$.  To my knowledge, the presence of the outer automorphism was
realized surprisingly recently, by van der Geer in 1982 (in terms of
the two isomorphisms of $Sp(4, \F_2)$ with $\SSS_6$, \cite[p.~323, 335, 337]{vdg},
see also \cite[p.~122]{do}):
$$
w_{\aaa}+\cdots + w_{\ff} = 0, \quad ( w_{\aaa}^2 + \cdots + w_{\ff}^2)^2 - 4 (w_{\aaa}^4 + \cdots + w_{\ff}^4) = 0.
$$
(Igusa's original equation \cite[p.~400]{igusa} obscured the
$\SSS_6$-action.)  Via the Gale transform (also known as the
``association map''), this is birational to the space of six points on
$\proj^1$, where six distinct points in $\proj^1$ induce six points on
$\proj^3$ by placing them on a rational normal curve (e.g.\ via $p_i
\mapsto [1;p_i; p_i^2; p_i^3]$), and six general points on $\proj^3$
induce six points on $\proj^1$ by finding the unique rational normal
curve passing through them.  
We describe the rational map $(\proj^1)^6 \dashrightarrow \Igusa$,
and then the rational map $(\proj^3)^6 \dashrightarrow \Igusa$.
\label{p3}\lremind{p3}

The rational map
$(\proj^1)^6 \dashrightarrow \Igusa$ is described as follows, using
the first description of the outer automorphism, \S
\ref{outer1}:
$$
W_{\xx} = \sum_{\{ A, \dots, E \} = \{ 1, \dots, 5 \}, \; \{ \al, \be,
  \ga \} = \{ 0,1,2 \}} N_{A, \dots, E} (p_6 p_A)^{\al}
( p_B p_C)^{\be} (p_D p_E)^{\ga}$$ where $N=2$ if the edge $BC$ has
the same color as edge $DE$, and $N=-1$ otherwise.  (A quick
inspection of the mystic pentagons shows that $\sum W_{\xx} = 0$.)
Hence the $\SSS_6$-representation on $H^0(\Igusa, \oh(1))$ is $O_5$,
the outer automorphism representation.  In
terms of the usual description of the outer automorphism \S \ref{usual},
$N=-1$ if $6A/BC/DE$ is a syntheme in the pentad, and $2$ otherwise.
\label{igusarocks}\lremind{igusarocks}

The birationality to the Segre cubic $\Segre$ arises by projective
duality ($\Segre$ and $\Igusa$ are dual hypersurfaces), which should
not involve the outer automorphism.  Indeed, the duality map $\Segre
\dashrightarrow \Igusa$ is given by
$$
W_{\xx} = Z_{\xx}^2 - \frac 1 6  \sum_{\yy=1}^6 Z_{\yy}^2,
$$
and the duality map $\Igusa \dashrightarrow \Segre$ is given by 
$$
Z_{\xx} = \left(  \sum_{\yy=1}^6 W_{\yy}^2 \right) W_{\xx} - 4 W_{\xx}^3 + \frac 2 3 \sum_{\yy=1}^6 W_{\yy}^3.
$$
It is perhaps surprising that these moduli maps are somehow ``dual'', while
the corresponding $\SSS_6$-representations $H^0(\Segre, \oh(1)) \cong O'_5$
and $H^0(\Igusa, \oh(1)) \cong O_5$ are {\em not} dual.  However,
their projectivizations are dual, as they differ by a sign representation
$\ep$.

The moduli map $(\proj^3)^6 \dashrightarrow \Igusa$ is frankly less
enlightening, but even here the outer automorphism perspective
simplifies the explicit formula.  Suppose the six points in $\proj^3$ are given by $[w_i; x_i;
y_i; z_i]$ ($1 \leq i \leq 6$). The usual invariants of this Geometric
Invariant Theory
quotient, in terms of tableaux, each have $9624$ monomials.  In terms
of the $Z$-variables, we have group orbits of $9$ monomials:
\begin{eqnarray*}
Z_{\xx} &=& \sum_{(\si, \tau) \in \SSS_5 \times \SSS_4} 
(\si, \tau) \circ \Big(
\frac 1 2 w_2 w_4 w_6 x_1 x_2 x_4 y_1 y_3 y_5 z_3 z_5 z_6
+ w_1 w_2 w_4 x_5 x_6^2 y_1 y_2 y_5 z_3^2 z_4  \\
& &  - \frac 1 2 w_2 w_3^2 x_5^2 x_6 y_2 y_4^2 z_1^2 z_6
+ 2 w_2 w_3 w_4 x_3 x_5 x_6 y_4 y_5 y_6 z_1^2 z_2
- w_1 w_2 w_4 x_3^2 x_4 y_5^2 y_6 z_1 z_2 z_6 \\
& &
- \frac 2 3 w_2 w_5 w_6 x_3^2 x_6 y_1^2 y_5 z_2 z_4^2
- \frac 1 2 w_1 w_2 w_3 x_1 x_5 x_6 y_2 y_3 y_4 z_4 z_5 z_6 \\
& & + \frac 1 6 w_2 w_3 w_4 x_1 x_2 x_5 y_1 y_4 y_6 z_3 z_5 z_6
+ \frac 1 4 w_1^2 w_2 x_2 x_3^2 y_5^2 y_6  z_4^2 z_6
\Big) .
\end{eqnarray*}
Here $\SSS_5$ acts by the ``outer action'' corresponding to $\xx$
on the six points $\{1, \dots, 6 \}$,
and $\SSS_4$ acts by permuting the co-ordinates $\{ w,x,y,z \}$
and by sign.
(There is significant abuse of notation in the way the formula
is presented, but hopefully the meaning is clear.)
This formula is less horrible than it appears, as the summands
can be interpreted as attractive geometric configurations on the
icosahedra of \S \ref{outer3}.

\bpoint{Six points in $\proj^2$} Finally, we describe the invariants
of six points in $\proj^2$ in terms of the mystic pentagons.  This
quotient is a double cover of $\proj^4$ branched over the Igusa
quartic $\Igusa$.  The Gale transform sends six points in $\proj^2$ to
six points in $\proj^2$, and exchanges the sheets.  The branch locus
of this double cover (the ``self-associated'' sextuples in the language of the Gale transform)
corresponds to when the six points lie on a conic;
by choosing an isomorphism of this conic with $\proj^1$, the rational
map $(\proj^1)^6 \dashrightarrow \Igusa$ is precisely the moduli map
described above in \S \ref{p3}.

Suppose the points in $\proj^2$ are $[x_i; y_i; z_i]$ ($1 \leq i \leq
6$).  We describe the moduli map $(\proj^2)^6 \dashrightarrow \proj^4$
in terms of the mystic pentagons \S \ref{outer1}:
$$W_{\xx} = \sum_{ \{ A, \dots, B \} = \{ 1, \dots, 6 \} }
N_{A, \dots, F} (x_A x_B) (y_C y_D) (z_E z_F).$$ Corresponding to each
term are two edges (corresponding to the pairs $AB$, $CD$, $EF$ not
containing $6$).  Then $N = 2$ if the two edges have the same color,
and $-1$ otherwise.  Notice the similarity to the moduli map for the
Igusa quartic above, in \S \ref{igusarocks}; this is not a coincidence,
and we have chosen the variable names $W_{\xx}$ for this reason.
Again, $N=-1$ if $6A/BC/DE$ is a syntheme in the pentad, and $2$
otherwise.

The condition for six points to be on a conic is for their image on
the Veronese embedding to be coplanar, hence that the following
expression is $0$:
$$
V := \det \left( \begin{array}{cccccc}
x_1^2 & y_1^2 & z_1^2 & x_1 y_1 & y_1 z_1 & z_1 x_1  \\
\vdots & \vdots &  \vdots & \vdots & \vdots & \vdots \\
x_6^2 & y_6^2 & z_6^2 & x_6 y_6 & y_6 z_6 & z_6 x_6  \\
\end{array}
\right).
$$
(Alternatively, the vanishing of this determinant ensures
the existence of a nontrivial quadric
$$
\al_{x^2} x^2 + \al_{y^2} y^2 + \al_{z^2} z^2 +
\al_{xy} xy + \al_{yz} yz + \al_{zx} zx = 0$$
satisfied by $(x_i, y_i, z_i)$ for all $i$ between $1$ and $6$.)

Then the formula for the double fourfold that is the
Geometric Invariant Theory quotient of six points on $\proj^2$ is 
$$\left( \sum W_{\xx}^2 \right)^2 - 4 \sum W_{\xx}^4 + 324 V^2 = 0$$
and it is clear that it is branched over the Igusa quartic $\Igusa$.
(See \cite[p.~17, Example~3]{do} for more information.)

\bpoint{Relation to the usual description of the invariants of six
  points on $\proj^1$} We relate the explicit invariant theory of \S
\ref{p1} to the classical or usual description of the invariants of
six points on $\proj^1$.  In the matching diagram language of
\cite{hmsv} (basically that of Kempe in 1894), the ring of
invariants is generated by the variables
$$
X_{\vec{AB} \cdot \vec{CD} \cdot \vec{EF} } =
(p_B-p_A) (p_D-p_C)(p_F-p_E)$$
where $\{ A, \dots, F \} = \{ 1, \dots, 6 \}$.
The variables $Z_{\mathbf{x}}$ of the Segre cubic threefold $\Segre$
are related to the matching diagrams in a straightforward way:
$$
X_{\vec{13} \cdot \vec{26} \cdot \vec{45} } =
(Z_{\mathbf{a}}+Z_{\mathbf{b}})/2$$
(and similarly after application of the
$\SSS_6$-action on both sides).  Notice that under the outer
automorphism, pairs are exchanged with ``synthemes'' (= partitions into
three pairs), and that is precisely what we see here.

This can of course be easily inverted, using
$$
Z_{\mathbf{a}} = ( Z_{\aaa} + Z_{\bb})/2 + (Z_{\aaa} + Z_{\cc})/2 -
(Z_{\bb} + Z_{\cc})/2.$$
As the $X$-variables form a $15$-dimensional vector space with
many relations, there are many formulas for the $Z$-variables in
terms of the $X$-variables.  

} 

\end{document}